\documentstyle{article}
\normalbaselines
\begin{document}
\bibliographystyle{unsrt}

\vbox{\vspace{12mm}}
\begin{center}
{\LARGE \bf ON AN ANALOG OF SELBERG'S }\\ [6mm]
{\LARGE \bf  EIGENVALUE CONJECTURE FOR $SL_{3}({\bf Z})$} \\ [10mm]
 


Sultan Catto$^{1}$\\{\it The Graduate School and University Center, and Baruch College\\ The City University Of New York \\ 17 Lexington Avenue \\ New York, NY 10010 \\ (e-mail: catto@gursey.baruch.cuny.edu)\\ and\\Physics Department, The Rockefeller University\\
1230 York Avenue, New York, NY 10021-6399\\(e-mail: 
cattos@rockvax.rockefeller.edu)}\\[8mm]

Jonathan Huntley$^{2}$ \\
{\it Department of Mathematics, Baruch College, CUNY \\
17 Lexington Avenue, New York, NY 10010\\
(e-mail: huntley@gursey.baruch.cuny.edu)}\\[8mm]

Jay Jorgenson$^{3}$\\
School of Mathematics\\
Institute for Advanced Study\\
Princeton, NJ 08540\\[8mm]

David Tepper$^{4}$\\
{\it Department of Mathematics, Baruch College, CUNY \\
17 Lexington Avenue, New York, NY 10010\\
(e-mail:tepper@gursey.baruch.cuny.edu)}\\[8mm]

\end{center}

\begin{flushleft}

(1) Work supported in part by DOE contracts No. DE-AC-0276-ER 03074 and 03075; 
NSF Grant No. DMS-8917754; and PSC-CUNY Research Award Nos. 6-6(6407, 7418, 8445).

(2) Work supported by PSC-CUNY Research Award No. 9203393 and NSF Award No. 
DMS-9200317.\\
(3) Work supported from NSF grant DMS-93-07023 and from the Sloan Foundation.\\
(4) Work partially supported by PSC-CUNY Research Award No. 62335.\\
\end{flushleft}
\newpage
\large{

\begin{abstract} 
\vskip .10in
Let $\cal H$ be the homogeneous space associated to the group
$PGL_3(\bf R)$. Let $X = \Gamma{\backslash {\cal H}}$ where $\Gamma =
SL_{3}(\bf Z)$ and consider the first non-trivial eigenvalue
$\lambda_1$ of the Laplacian on $L^2(X)$.  Using geometric
considerations, we prove the inequality $\lambda_1 > 3\pi^2/10$.
Since the continuous spectrum is represented by the band
$[1,\infty)$, our bound on $\lambda_{1}$ can be viewed as an analogue
of Selberg's eigenvalue conjecture for quotients of the hyperbolic
half space.   

\vskip .10in
\noindent
\end{abstract}
\vskip 1cm
{\bf \S 1. Statement of the main theorem}
\vskip .10in

A fundamental question in the spectral theory of automorphic forms is
whether small eigenvalues exist.  More specifically, let $G$ be a
noncompact reductive group with finite center, $\Gamma$ a nonuniform
lattice, $K$ a maximal compact subgroup of $G$, and set $X =
\Gamma\backslash G/K$.  It is well known from the theory of
Eisenstein series that $L^2(X)$ has continuous spectrum for the
ring of invariant differential operators, and in particular for the
positive Laplacian, $\bf \Delta$.  The continuous spectrum will be,
in cases of interest such as $PGL_n(\bf R)$, an interval $[a,\infty)$
with $a > 0$.  The question we referred to above is:  Do non-constant
square integrable eigenforms exist with eigenvalue $\lambda < a$?
This problem is important for various considerations in number
theory.  In the case $G=PGL_{2}({\bf R})$ and $\Gamma$ is a
congruence subgroup, Selberg conjectured that no such nontrivial small
eigenvalues exist.  

\vskip .10in
In this paper, we consider the case when 
$$
G = PGL_3({\bf R})~~~~ and~~~~\Gamma
= SL_3({\bf Z}).
$$
Our main result is the following.
\vskip .10in
{\,\,\,\,\,\bf Theorem} {\it
Let $\lambda_1$ be the eigenvalue for the first nontrivial
eigenform on $L^2(X)$.  Then} 
$$
\lambda_1 > 3\pi^{2}/10 > 2.96088.
$$

\vskip .10in
\noindent
{\bf \S 2. Notation.}
\vskip .10in

Let $\cal H = G/K$ and set $X = \Gamma\backslash \cal H$.  Explicit
coordinates for $\tau \in \cal H$ via the Iwasawa decomposition are given 
by  
\[
\tau = 
\left(\begin{array}{ccc}
 1&x_{2} & x_{3}  \\ 0 & 1 & x_1 \\ 0 & 0 &1
\end{array}\right)
\left(\begin{array}{ccc}
 y_1 y_2 & 0 & 0 \\ 0 & y_1 & 0 \\ 0 & 0 &1
\end{array}\right)
\]    
with $y_{1}, y_{2} > 0$, 
from which one can compute that the (positive) Laplacian $\bf \Delta$
can be written as
$$
-{\bf \Delta} = y^2_1 \frac{\partial^2}{\partial y^{2}_{1}} -
y_1y_2\frac{\partial^2}{\partial y_{1}\partial y_{2}} + y^2_2
\frac{\partial^2}{\partial y_{2}^{2}} 
+ y^2_1(x^2_2 + y^2_2) \frac{\partial^{2}}{\partial x^{2}_{3}} +
y^2_1 \frac{\partial^{2}}{\partial x_{1}^{2}} + y^2_2
\frac{\partial^{2}}{\partial x^{2}_{2}}
+ 2y^2_1 x_2 \frac{\partial^{2}}{\partial x_{1}\partial
x_{3}}
$$
(see pages 17 and 33 of ${\cite{Bu 84}}$).  The ring of invariant
differential operators is spanned by the identity, the Laplacian
$\bf \Delta$, and a third order operator $\bf \Delta_3$ (see ${\cite{Bu 84}}$).  The invariant volume element is given 
by 
$$
dV =\frac{dx_1dx_2dx_3dy_1dy_2}{y^{3}_{1}y^{3}_{2}}.
$$ 
We shall not use the explicit formula for the invariant volume
element; however, the above expression for the Laplacian will be
necessary in our proof of the main theorem.

\vskip .10in
For our purposes, it is more convenient to work with functions on
${\cal H}$ that are $SL_{3}({\bf Z})$ invariant rather than considering
functions on the quotient space $X$.  To this end, we introduce a
fundamental domain for $\Gamma\backslash \cal H$. Specifically,
computations on page 56 of ${\cite{Gr 93}}$ show that a fundamental
domain $\cal D$ is described through the following set of
inequalities: 

\qquad \qquad
\begin{tabular}{l}
$v^{\frac{3}{2}} < v^{\frac{3}{2}} (1-x_2 + x_3)^2 + w(1-x_1)^2 +w^{-1};$  \\
$v^{\frac{3}{2}} < v^{\frac{3}{2}}(x_2-x_3)^2 + w(1-x_1)^2 + w^{-1}_1;$ \\
$v^{\frac{3}{2}} < v^{\frac{3}{2}} x^2_2 + w; 
\,\,\,\,\,\,\,
v^{\frac{3}{2}} < v^{\frac{3}{2}} x^2_3 + wx^2_1 + w^{-1};$ \\
$1 < w^{-2} + x^2_1; \,\,\,\,\,\,\,
0 < x_1 < \frac{1}{2}; \,\,\,\,\,\,\,
0 < x_2 < \frac{1}{2} \,\,\,\,\,\,\,
-\frac{1}{2} < x_3 < \frac{1}{2},$ \\ 
\end{tabular} 
\flushleft
where we have used the notation $w^{-1} = y_1$ and $v^{-\frac{3}{2}}
= y^2_2 y_1$.  Let $S$ denote the Siegel set described via the
inequalities 
$$
0 < x_1 < \frac{1}{2}, \,\,\,\,\,
0 < x_2 < \frac{1}{2}, \,\,\,\,\,
-\frac{1}{2} < x_3 < \frac{1}{2}, \,\,\,\,\,
y_1 > \frac{\sqrt{3}}{2},\,\,\,\,\, 
y_2 >\frac{\sqrt{3}}{2}. 
$$
The set $S$ contains the fundamental domain $\cal D$.  Further,
results from page 61 ${\cite{Gr 93}}$ show the existence of elements
$\gamma_{1}, ..., \gamma_{10} \in SL_{3}({\bf Z})$ such that
$S \subset \bigcup\limits^{10}_{i=1} {\cal D} {\gamma_{i}}$ (we have used
the notation $\cal D \gamma$ to denote the image of the fundamental
domain $\cal D$ under left multiplication by $\gamma$).
The main aspects of the above points which we shall use are the
assertions that for any $\tau \in S$ we have $y_{1}(\tau) >
\sqrt{3}/2$ and that $S$ is contained in ten translates of $\cal D$.

\vskip .10in
Recall that an automorphic form is a $C^{\infty}$ function $\phi$ on
$\cal H$ which satisfies the following properties:
\begin{quote}
(i)~ $\phi (\gamma\circ\tau) = \phi(\tau)$ for $\gamma\in
SL_3(\bf Z)$;

(ii)~ $\vert\phi (\tau)\vert \ll y^{N_{1}}_1 y^{N_{2}}_2$
for $\tau\in \cal D$ and integers $N_1, N_2$;

(iii)~ $\phi$ is an eigenform for the ring of invariant
differential operators.
\end{quote}
\noindent
An automorphic form is said to be a cusp form if it satisfies the
additional property
\begin{quote}
(iv)~
$\displaystyle
\int\limits^{\frac{1}{2}}_{-\frac{1}{2}}
\int\limits^{\frac{1}{2}}_{-\frac{1}{2}} \phi\left(
\left[\begin{array}{ccc}
 1 & 0 &\xi_3 \\ 0 & 1 & \xi_1\\ 0 & 0 & 1
\end{array}\right]
 \tau\right) d\xi_1
d\xi_3 
=\int\limits^{\frac{1}{2}}_{-\frac{1}{2}}
\int\limits^{\frac{1}{2}}_{-\frac{1}{2}}  \phi\left(
\left[\begin{array}{ccc}
 1 & \xi_2 & \xi_3 \\ 0 & 1 & 0 \\ 0 & 0 & 1
\end{array}\right]
\tau\right) d\xi_2
d\xi_3 = 0 $
\end{quote}
\noindent
Cusp forms are square integrable.  Although we shall not need this
fact, let us note that, from the theory of Eisenstein series, the
only noncuspidal square integrable automorphic forms on $X$ are
constant.   

\vskip .10in
\noindent
{\bf \S 3. Proof of the main theorem}
\vskip .10in
Our method of proof is a modification of that used by Roelcke to
show that the small eigenvalue $\lambda_{1}$ for the quotient space
$SL_2({\bf Z})\backslash SL_2({\bf R})/ SO_2({\bf R})$
satisfies the bound $\lambda_1 > 3 \pi^{2}/2$ (see page 511 of
${\cite{He 83}}$).  We shall use the Fourier expansion of automorphic
forms associated to $SL_{3}(\bf Z)$, as developed in Chapter IV of
${\cite{Bu 84}}$.  

\vskip .10in
Assume that $\phi$ is a non-constant automorphic form, so then
${\bf \Delta}\phi =\lambda\phi$ and $\phi{\bf \Delta}\phi =
\lambda\phi^2$.   Through integration by parts, using the automorphic
boundary conditions, and the fact that the Siegel domain $S$ is
contained in ten translates of the fundamental domain $\cal D$, we
obtain the inequality  
$$
\frac{\int\limits_S\vert\nabla\phi\vert^{2} dV}
{\int\limits_S \vert\phi\vert^{2} dV}  < 10 \lambda.
$$
As on page 67 of ${\cite{Bu 84}}$, let us expand $\phi$ in a Fourier
expansion with respect to the abelian group
$$
\left\{
\left[
\begin{array}{ccc}
 1 & 0 & \xi_3 \\ 0 & 1 & \xi_1 \\ 0 & 0 &1
\end{array}
\right]
 \,\,\,\,\,\xi_{1}, \xi_{3} \in \bf R \right\}.
$$
Specifically, we have $\displaystyle\phi (\tau ) =
\sum\limits_{n_{1},n_{3}}\phi^{n_{3}}_{n_{1}} (\tau)$ 
where
$$
\phi^{n_{3}}_{n_{1}}(\tau ) = \int\limits^1_0 \int\limits^1_0
\phi\left(
\left[
\begin{array}{ccc} 
1 & 0 & \xi_3\\ 0 & 1 & \xi_1 \\ 0 & 0 & 1
\end{array}
\right)
\tau\right] e^{-2\pi i(n_{1}\xi_{1}+n_{3}\xi_{3})} d\xi_1
d\xi_3.
$$
Observe that $\phi^{0}_{0}=0$ since $\phi$ is not constant and square
integrable, hence cuspidal.  Let
$$
\Gamma^{2}_{1} = \left\{
\left[
\begin{array}{ccc}
 a & b & 0 \\ c & d & 0 \\ 0 & 0 &1
\end{array}
\right]
 : 
\left[
\begin{array}{cc}
 a & b \\ c & d
\end{array}
\right]
 \in
SL_2({\bf Z})\right\},
$$
and set $\Gamma^2_\infty$ to be the subgroup of $\Gamma^{2}_{1}$ which
stablizes infinity.  As on page 69 of ${\cite{Bu 84}}$, we then can write
$\phi^{n_{3}}_{n_{1}}$ as 
$$
\phi^{n_{3}}_{n_{1}} (\tau) =
\sum\limits_{\gamma\in\Gamma^{2}_{\infty}\backslash\Gamma^{2}_{1}}
\sum\limits_{n_{1}=1}^{\infty}\phi^0_{n_{1}}(\gamma\circ\tau).
$$  
By a standard application of elliptic regularity, $\phi$ is
necessarily $C^\infty$, hence we can interchange integration and
summation and apply Parseval's theorem to obtain the inequality
$$
\frac{\int\limits_{S}
\sum\limits_{n_{1}=1}^{\infty}\left|
\sum\limits_{\gamma\in\Gamma^{2}_{\infty}\backslash\Gamma^{2}_{1}}
\nabla_{\tau}\phi^0_{n_{1}}
(\gamma\circ\tau)\right|^2 dV}
{\int\limits_{S}
\sum\limits_{n_{1}=1}^{\infty}
\left|
\sum\limits_{\gamma\in\Gamma^{2}_{\infty}\backslash\Gamma^{2}_{1}}
\phi^0_{n_{1}}(\gamma\circ\tau)\right|^2 dV} < 10 \lambda.
$$
Since $\nabla$ is an invariant operator, we may differentiate the
expressions in the numerator and then evaluate at $\gamma\circ\tau$,
thus yielding
$$
\frac{\int\limits_{S}
\sum\limits^\infty_{n_{1}=1}\left|
\sum\limits_{\gamma\in\Gamma^{2}_{\infty}\backslash\Gamma^{2}_{1}}
\nabla_{\tau}\phi^0_{n_{1}}(\tau)\big|_
{\gamma\circ\tau}\right|^2 dV}
{\int\limits_S\sum\limits^\infty_{n=1}
\left|\sum\limits_{\gamma\in\Gamma^{2}_{\infty}\backslash\Gamma^{2}_{1}}
\phi^0_{n_{1}}(\tau)\big|_{\gamma\circ\tau}\right|^2 dV} < 10 \lambda.
$$

We now integrate by parts and consider the action of the Laplacian
$\bf \Delta$ on functions of the form $\phi_{n_{1}}^{0}$.  Since each
function $\phi_{n_{1}}^{0}$ is independent of $x_{3}$, these terms in
$\bf \Delta$ annihilate $\phi_{n_{1}}^{0}$.  Observe that all terms
involving $y_{1}, y_{2}$ and $x_{2}$ are positive operators (compare
with line (2.31) on page 32 of ${\cite{Bu 84}}$), so we obtain the bound
$$
{\bf \Delta} \phi_{n_{1}}^{0} \geq - y_{1}^{2} \cdot
\frac{\partial^{2}}{\partial x_{1}^{2}}\phi_{n_{1}}^{0} 
= y_{1}^{2} \cdot 4 \pi ^{2}n_{1}^{2}\phi_{n_{1}}^{0}.
$$
Since $y^2_1 > \frac{3}{4}$, we have
${\bf \Delta} \phi_{n_{1}}^{0} \geq \frac{3}{4}\cdot 
4 \pi ^{2}n_{1}^{2}\phi_{n_{1}}^{0}=
3 \pi ^{2}n_{1}^{2}\phi_{n_{1}}^{0}$.  Combining this inequality with
the above calculations and the cuspidality condition $\phi_{0}^{0} = 0$,
we obtain 
$$
10 \lambda >  
\frac{ \int\limits_{S}
\sum\limits^\infty_{n_{1}=1}
\sum\limits_{\gamma\in\Gamma^{2}_{\infty}\backslash\Gamma^{2}_{1}}
{\bf \Delta}_{\tau}\phi^0_{n_{1}}(\tau)\big|_{\gamma\circ\tau} 
\cdot \phi^0_{n_{1}}(\tau)\big|_{\gamma\circ\tau} dV}
{\int\limits_S\sum\limits^\infty_{n=1}
\left|\sum\limits_{\gamma\in\Gamma^{2}_{\infty}\backslash\Gamma^{2}_{1}}
\phi^0_{n_{1}}(\tau)\big|_{\gamma\circ\tau}\right|^2 dV} \geq 3\pi^{2},
$$
hence $\lambda\geq 3\pi^2/10$.  Since $\phi$ was any cusp form, we
obtain the bound as asserted in the theorem.

\vskip .10in
{\bf \S 4. Concluding remarks}
\noindent
\vskip .10in

As the continuous spectrum in this situation is $[1,\infty)$, our
theorem implies an analogue of Selberg's eigenvalue conjecture.
Note that our bound is stronger than result for $SL_{3}(\bf Z)$
from ${\cite{Mi 96}}$ who proved $\lambda_{1} \geq 1$.  In general, our
method applies to $G = SL_{n}(\bf R)$ with $\Gamma = SL_{n}(\bf Z)$ to give the bound $\lambda_{1} > 3\pi^{2}/M$ where $M$ is the
number of fundamental domains which intersect a Siegel set containing
the fundamental domain constructed in ${\cite{Gr 93}}$; however, for $n
\geq 4$, this bound is rather weak.  Finally, let us remark that our
theorem is indeed a consequence of the Ramanujan conjecture, which
asserts that all nontrivial automorphic representations come from
tempered representations.    

\noindent

\vfill

}
\enddocument